\documentclass{elsart}
\usepackage{times, mathptm}  
\usepackage{amsfonts, epsfig}

\newcommand{\C}{{\Cset}}
\renewcommand{\H}{{\Hset}}
\newcommand{\Z}{{\Zset}}

\newcommand{\PSL}{ {\mathrm{PSL}} }
\newcommand{\SL}{ {\mathrm{SL}} }

\newcommand{\bdry}{\partial}
\newcommand{\cross}{\times}
\newcommand{\from}{{\colon\thinspace}}
\newcommand{\injects}{{\hookrightarrow}}

\newcommand{\tr}{ {\mathrm{\, tr}} }
\newcommand{\SnapPea}{SnapPea}

\renewcommand{\a}{\alpha}
\renewcommand{\b}{\beta}
\renewcommand{\t}{\tau}
\newcommand{\g}{\gamma}
\newcommand{\mtext}[1]{\quad\mbox{#1}\quad}

\newcommand{\CP}{{ {\mathbb C}\, 
    \mbox{\fontfamily{cmr}\fontshape{n}\selectfont P} }}

\renewcommand{\phi}{\varphi}

\begin{document}

\begin{frontmatter}
\title{Examples of non-trivial roots of unity at ideal points of
  hyperbolic 3-manifolds }

\author{ Nathan M.~Dunfield\thanksref{NSF} }
\address{ Department of Mathematics, University of Chicago,
  5734 S.~University Ave., Chicago, IL 60637, USA,
 e-mail: nathand@math.uchicago.edu}

\thanks[NSF]{This work was partially supported by an NSF Graduate
  Fellowship.}

\begin{abstract}
 
This paper gives examples of hyperbolic 3-manifolds whose $\SL(2,\C)$
character varieties have ideal points for which the associated roots
of unity are not $\pm 1$.  This answers a question of Cooper, Culler,
Gillet, Long, and Shalen as to whether roots of unity other than $\pm
1$ occur.
 
\end{abstract}
\end{frontmatter} 

\section{Introduction}

Constructing incompressible surfaces in 3-manifolds via degeneration
of hyperbolic structures was introduced by Culler and Shalen in
\cite{CullerShalen83}.  This technique has been very useful.  For
instance, it was the basis for half of the proof of the Cyclic Surgery
Theorem \cite{CGLS}.  The construction starts with a complex curve of
representations of the fundamental group of the manifold into
$\SL(2,\C)$.  The surfaces are constructed from ideal points of this
curve.  In \cite{CCGLS}, Cooper, Culler, Gillet, Long, and Shalen
showed that certain ideal points have roots of unity associated to
them which reflect the topology of the associated surfaces.  They
asked whether these roots are always $\pm 1$, this being the case in
all known examples.  I will give examples of hyperbolic manifolds
where some of these roots are not $\pm 1$, but are fourth or sixth
roots.

In Section~\ref{background}, I review the construction of
incompressible surfaces from ideal points, and the associated roots of
unity.  I found the examples using a heuristic method described in
Section~\ref{motivation}.  I discuss examples which appear to have
ideal points where the root of unity is a fourth or sixth root in
Section~\ref{probable_examples}.  I will rigorously show that two of
these probable examples have ideal points whose associated roots of
unity are not $\pm 1$.  The first is a punctured torus bundle where
the root of unity is $\pm i$.  The second example is the complement of
a knot in $S^2 \cross S^1$ and the root of unity is a sixth root.
These examples are treated in Sections~\ref{N} and \ref{M},
respectively.  The first example, $N$, also has two other interesting
properties which are discussed in the remarks at the end of
Section~\ref{N}.  In the language of Section~\ref{background}, let
$X(N)$ be the character variety of $N$ and $X_0$ a component of $X(N)$
which contains the character of a discrete faithful representation.
The manifold $N$ has a strict boundary slope which is not associated
to any ideal point of $X_0$ but which is associated to an ideal point
of some other component of $X(N)$.  Also, there is a $\gamma$ in
$\pi_1(N)$ such that $\tr_\gamma$ is constant on $X_0$. 

I would like to thank Peter Shalen for suggesting this problem to me
and for innumerable useful discussions.  I would also like to thank
the referee for very useful comments and suggestions.   I should also
mention that, in unpublished work, Cooper and Long found a hyperbolic
manifold with an ideal point where the root of unity $\pm i$.

\section{Background}\label{background}

The basic reference is \cite{CGLS}.  Let $M$ be a compact 3-manifold.
Let $X$ denote its $\SL(2,\C)$ character variety.  Basically, $X$ is
the set of representations of $\pi_1(M)$ into $\SL(2,\C)$ modulo
conjugacy.  It has a natural structure as an affine algebraic variety
over $\C$.  Let $X_0$ be a curve in $X$.  For instance, if $M$
hyperbolic with one cusp, an irreducible component of $X$ which
contains the character of a discrete faithful representation is such a
curve.  Let $f\from \widetilde X_0 \to X_0$ be a birational map with
$\widetilde X_0$ a smooth projective curve.  The points in $\widetilde
X_0 \setminus f^{-1}(X_0)$ are called ideal points, and to each ideal
point there is an associated action of $\pi_1(M)$ on a simplicial
tree.  This action yields an incompressible surface by choosing a
suitable equivariant map from the universal cover of $M$ to the tree,
$T$, taking the inverse image of the midpoints of the edges $T$, and
pushing the resulting surface down to $M$.

Now suppose that the boundary of $M$ is a torus.  A slope is an
isotopy class of simple closed curves in $\bdry M$, or equivalently, a
pair of primitive elements $\{ \alpha, -\alpha \}$ in $\pi_1(\bdry
M)$.  If a surface properly embedded in $M$ has non-empty boundary,
then its boundary consists of disjoint simple closed curves in $\bdry
M$.  Since these curves are parallel, they have the same slope, which
is called the boundary slope of the surface.  For each $\gamma \in
\pi_1(M)$ there is a natural regular function $\tr_\gamma$ on $X$
whose value at a character of a representation $\rho$ is $\tr
(\rho(\gamma))$.  These induce meromorphic functions on $\widetilde
X_0$.  Let $x$ be an ideal point.  If the trace function of some
peripheral element has a pole at $x$, then there is exactly one slope
$\{ \alpha, -\alpha \}$ for which $\tr_{\pm \alpha}(x)$ is finite.
Any incompressible surface $S$ associated to $x$ has non-empty
boundary with boundary slope $\{\alpha, -\alpha \}$.  Suppose the
number of boundary components of $S$ is minimal among all
incompressible surfaces associated to $x$.  For such $S$, it was shown
in \cite{CCGLS} (for another proof see \cite{CooperLong93}) that
$\tr_\alpha(x) = \lambda + \lambda^{-1}$ where $\lambda$ is a root of
unity whose order divides the number of boundary components of any
connected component of $S$.  I will call $\lambda$ the root of unity
associated to $x$.

\section{Heuristic method}\label{motivation}

Before explaining the method, let me summarize the
connection between the character variety and the Thurston deformation
variety $T$ (see \cite{ThurstonLectureNotes} and \cite{NeumannZagier},
or \cite{BenedettiPetronio} for details).  Thurston defined $T$ in his
proof of the Hyperbolic Dehn Surgery Theorem.  Fix a triangulation of
$M \setminus \bdry M$ by ideal tetrahedra.  I mean this topologically,
i.e.~$M \setminus \bdry M$ decomposed as the union of 3-simplices
minus their vertices -- these tetrahedra should not yet be thought of
as sitting in hyperbolic three space $\H^3$.  Up to isometry, an ideal
tetrahedron in $\H^3$ with geodesic sides and a preferred edge is
parameterized by a single complex number, the edge parameter, which
determines the similarity type of the triangle which is the link of an
ideal vertex.  If we choose an edge parameter for each tetrahedron in
our triangulation, there are algebraic conditions, the gluing
equations, which hold if the hyperbolic structures on the tetrahedra
glue up to give a hyperbolic structure on $M$.  The space of edge
parameters that satisfy the gluing equations is the variety $T$.  A
point $p$ in $T$ where none of the edge parameters is totally
degenerate ($= 0, 1,$ or $\infty$) defines a representation of
$\pi_1(M)$ into $\PSL(2, \C)$ by means of a simplicial ``developing''
map. If all the imaginary parts of the edge parameters have the same
sign, and the gluing equations still hold in the universal cover of
$\C \setminus \{0 \}$, then $p$ defines a hyperbolic structure on $M$.
In this case, the representation associated to $p$ is just the
holonomy representation of this structure.

There is a map from $T$ to the $\PSL(2,\C)$ (rather than $\SL(2, \C)$)
character variety, by sending a point to its associated
representation.  So there is almost a map from $T$ to $X$.  It is well
known that $T$ and $X$ are very similar.  In fact, at least on a
component of $T$ containing a complete structure, it is possible to
construct a finite-to-finite correspondence (in the sense of algebraic
geometry \cite{Mumford}) from $T$ to $X$, though I will not use this
explicitly.

For a point $p$ in $T$, in addition to the associated representation
$\rho\from \pi_1(M) \to \PSL(2, \C)$ there is a representation $h\from
\pi_1(\bdry M) \to \C$ coming from the similarity
structure of $\bdry M$ induced by the hyperbolic structure defined by
$p$.  Fixing $\alpha \in \pi_1(\bdry M)$ we get a rational function
$h_\alpha$ on $T$.  If we conjugate $\rho$ so that $\rho(\pi_1(\bdry M))$
fixes $\infty$ in $\bdry\H^3$,
\[ \rho(\alpha) = \left( \begin{array}{cc} \lambda & \delta \\
                                                                        0       & \lambda^{-1} 
\end{array} \right). \]
Then $h_\alpha(p)$ is $\lambda^2$ and we have:
 \[
 \tr(\rho(\alpha))= \pm \left( \sqrt{h_\alpha(p)} +
\left(\sqrt{h_\alpha(p)}\right)^{-1} \right).
\]

The idea for finding an example is this: Suppose we want to do 1/3
orbifold hyperbolic Dehn filling along a peripheral curve $\alpha$.
Recall from the proof of Thurston's Hyperbolic Dehn Surgery Theorem
that we do this by moving around in $T$ until we make $h_\alpha = e^{2\pi
  i/3}$.  If we can't do 1/3 orbifold hyperbolic Dehn filling, this
would mean that $h_\alpha$ is never $e^{2\pi i/3}$ for any point of
$T$ corresponding to a non-degenerate hyperbolic structure.  Let's
suppose this is the case.  Now $h_\alpha$ is a non-constant rational
function on $T$, and so it takes on $e^{2\pi i/3}$ at some point if we
replace $T$ by a projective closure.  Our assumptions make it likely
that this point corresponds to an ideal point.  I used Weeks' wonderful
program \SnapPea\ \cite{SnapPea} to search for manifolds where 1/3
filling was not possible.  Of the several thousand manifolds tested,
only ten had such a slope.  In nine of these cases, numerical evidence
suggested that as the program tried to do 1/3 filling, it was actually
getting sucked toward an ideal point where the holonomy was $e^{2\pi
  i/3}$, and therefore the associated root of unity was a sixth root.
There is nothing special about $1/3$ -- you can do the same thing for
any $1/n$.  For $n=2$, I found probable examples including the
meridian slope of the knot complements $8_{16}, 8_{17}, 9_{29}$, and
$9_{32}$, and a certain punctured torus bundle.

I will rigorously show that two of these examples have ideal points
whose associated roots of unity are not $\pm 1$.  For the punctured
torus bundle in Section~\ref{N}, it is possible to completely compute
$X$ by hand, and, while I found this example using the above heuristic
method, no use is made of $T$.  For the knot complement in $S^2 \cross
S^1$ I worked directly with $T$ and some computer computations are
needed.

\section{Probable examples} \label{probable_examples}

By searching the Callahan-Hildebrand-Weeks census of non-compact
orientable hyperbolic 3-manifolds that can be triangulated with seven
or fewer ideal tetrahedra (\cite{HildebrandWeeks} and
\cite{CallahanHildebrandWeeks}), as discussed above, I found the
following examples which probably have ideal points with non-trivial
sixth roots: m137, s783, s784, s938, s939, v2946, v2947, v3219, and
v3221.  Examples which appear to have ideal points where the root of
unity is $\pm i$ include: m135, m147, and the knot complements
$8_{16}, 8_{17}, 9_{29}$, and $9_{32}$.  In all cases, the boundary
slope of a surface associated with the ideal point is a slope that is
the shortest in the Euclidean structure of the cusp.  This is not
surprising given the Gromov-Thurston $2\pi$ theorem
\cite{BleilerHodgson96} and the fact that one can't do $1/3$ or $1/2$
orbifold Dehn filling.

An extra difficulty arrises when using the heuristic method to find
examples for which the root of unity is $\pm i$, rather than a sixth
root.  The problem is that it often happens that when the 1/2 filling
is not hyperbolic, there is still no ideal point with root of unity
$\pm i$.  For instance, one can't do 1/2 filling along the meridian of
the complement of a Montesinos knot because the 2-fold branched cover
of $S^3$ along the knot is Seifert fibered.  But if the meridian is to
be the boundary slope of an incompressible surface, the knot
complement must contain a closed incompressible surface \cite{CGLS}.
Montesinos knots have closed incompressible surfaces if and only if
they have more than three rational tangles \cite{Oertel84}.

\section{Punctured torus bundle}\label{N}

The first  example is the punctured torus bundle over $S^1$ with monodromy
\[\phi = \left( \begin{array}{cc} -1 & -2 \\ -2 & -5 \end{array}\right) =
        -\left(\begin{array}{cc} 1 & 0 \\ 1 & 1 \end{array}\right)^2
        \left(\begin{array}{cc} 1 & 1 \\ 0 & 1 \end{array}\right)^2,
\] 
which I will call $N$.  Since $\phi$ is hyperbolic, work of Thurston
shows that N has a complete hyperbolic metric of finite volume (see
\cite{Otal96}).  Let $X$ be the character variety of $N$, and let
$X_0$ be an irreducible component of $X$ which contains the character
of a discrete faithful representation.  It turns out that $X_0$ is
birationally isomorphic to $\CP^1$, has exactly four ideal points, and
the associated roots of unity are all $\pm i$.  Incompressible
surfaces in punctured torus bundles have been classified
(\cite{FloydHatcher82} and \cite{CullerJacoRubinstein}), and $N$
contains three besides the fiber.  Information about these surfaces is
given below; boundary slopes are with respect a basis where the fiber
has slope $\infty$.
\begin{center}
\begin{tabular}{|c|c|c|c|}
\hline
Name & Slope & \# of Boundary Comp. & Genus \\
\hline
$S_1$ & 0 & 4 & 0 \\
$S_2$ & 1 & 4 & 0 \\
$S_3$ & 1/2 & 2 & 1 \\
\hline
\end{tabular}
\end{center}
The surfaces associated to the ideal points of $X_0$ are $S_1$ and
$S_2$.  Thus $N$ is one of the few examples known where there is a
component of the character variety with only two surfaces associated
to all its ideal points (Culler and Shalen have developed a theory of
such manifolds \cite{CullerShalenTwoSurface}).  The Dehn filling of
$N$ along the slope 1/2 , $N_{1/2}$, is a manifold whose character
variety has dimension 1.  Ideal points of the character variety of
$N_{1/2}$ have associated surface $S_3$.  The slopes of $S_1$ and
$S_2$ are a basis for $\pi_1( \bdry N)$ and one can compute from $X_0$
that the cusp shape with respect to this basis is $i$.

I will use the methods of \cite[Chap.~6]{Hodgson86} and \cite[Sec.~4.5]{Porti97}.  See
\cite{CullerShalen83} for a precise definition of the character
variety.  Let $F$ be a once-punctured torus.  Choose generators $a$
and $b$ for the free group $\pi_1(F)$ so that the action of $\phi$ is
\[
\phi(a) = b^{-1} a^{-1} b^{-1} \mtext{and} 
\phi(b) = b a \left(b^{-1} a^{-1} b^{-1} \right)^3,
\] 
and so
\[
\pi_1(N) = \left\langle a, b, t \ \Big| \ t a t^{-1} = b^{-1} a^{-1}
b^{-1}, t b t^{-1} = b a \left( b^{-1} a^{-1} b^{-1} \right)^3\right\rangle.
\]
This presentation was chosen so that $t$ and $l = a b a^{-1} b^{-1}$
commute and are a basis for $\pi_1(\bdry N)$.  

An irreducible representation $\rho$ of $\pi_1(F)$ into $\SL(2,\C)$
 with $a \mapsto A$ and $b \mapsto B$ can be conjugated so that 
\[
A = \left( \begin{array}{cc} x & 1 \\ 0 & 1/x \end{array} \right)
\mtext{and} 
B = \left( \begin{array}{cc} y & 0 \\ z &1/y \end{array} \right),
\]
where $x$ and $y$ are arbitrary eigenvalues of $A$ and $B$
respectively.  Since $z$ is determined by $\tr(A B)$, the conjugacy
class of an irreducible representation of $\pi_1(F)$ is completely
determined by $\tr(A)$, $\tr(B)$, and $\tr(AB)$.  So, for such an
irreducible representation of $\pi_1(F)$ to extend to all of
$\pi_1(N)$ it is necessary and sufficient that
\begin{eqnarray}
\tr(A) &=& \tr(\rho(\phi(a)))  =  \tr( B^{-1} A^{-1} B^{-1}),  \label{relation1}\\
\tr(B) &=&  \tr( \rho(\phi(b)))  =\tr( B A (B^{-1} A^{-1} B^{-1} )^3), \label{relation2} \\
\tr(A B) &=&  \tr(\rho(\phi(ab))) = \tr(B^{-1} A^{-1} B^{-1} B A (B^{-1} A^{-1} B^{-1} )^3), \label{relation3}
\end{eqnarray}
as then there exists a $T$ in $\SL(2,\C)$ such that for all $g \in
\pi_1(F)$, we have $\rho( \phi(g) ) = T \rho(g) T^{-1}$.  We can then
take $\rho(t) = T$.  Note that $T$ is unique up to sign because the
stabilizer under conjugation of an irreducible representation is $\{
\pm I \}$.  

A representation $\rho$ of $\pi_1(F)$ is reducible if and only if
$\tr(\rho( l) ) = 2$.  I will only be interested in components of $X$
which contain the character of a discrete faithful representation.  By
Proposition 2 of \cite{CullerShalen84}, the function $\tr_l$ is
non-constant on such components.  Thus I will restrict attention to
the subvariety $X'$ of $X$ which is the Zariski closure of the subset
of characters which are irreducible when restricted to $\pi_1(F)$.

By Proposition~1.4.1 of \cite{CullerShalen83}, coordinates on $X'
\subset \C^7$ are $\a = \tr(A)$, $\b = \tr(B)$, $\g = \tr(A B)$, $\t =
\tr(T)$, $\tr(A T)$, $\tr(B T)$, and $\tr(A B T)$.      
The basic tool will be the following identities, which hold for all
$X, Y \in \SL(2, \C)$:
\begin{equation} \label{idents}
\tr(X Y) = \tr( Y X), \tr(X) = \tr(X^{-1}), \mtext{and} 
        \tr(X Y) = \tr(X) \tr(Y) - \tr(X Y^{-1}).
\end{equation}
Using these to expand relation (\ref{relation1}) gives
\[
\tr(A) = \tr(B^{-1} A^{-1} B^{-1}) =  \tr(B A B) = \tr(B) \tr(A B) - \tr(A)
\mtext{or}  2 \a = \b \g.
\] 
Similarly expanding (\ref{relation2}) and using (\ref{relation1}) to
note $\tr(A B^2) = \tr(A)$, we get $2 \b = \a ( \a \b - \g)$.  If
(\ref{relation1}) and (\ref{relation2}) hold then so does
(\ref{relation3}) because in this case
\begin{eqnarray*}
\tr(\rho(\phi(a b))) &=& \tr(\rho(\phi(a))) \tr(\rho(\phi(b))) - 
                \tr\left(\rho(\phi(a)) \rho(\phi(b))^{-1}\right)\\
&=& \tr(A) \tr(B) - \tr\left(B^{-1} A^{-1} B^{-1} \cdot B^{-1}\right) \\
&=& \tr(A)\tr(B) - \left(\tr\left(B^{-1} A^{-1} B^{-1}\right) \tr(B) - \tr(AB)\right) =
\tr(AB).
\end{eqnarray*}
Combining the information from (\ref{relation1}) and
(\ref{relation2}), we get that either $p \equiv \b^2 \g^2 - 2 \g^2 - 8
= 0$ or $\b = \a = 0$.  If $\rho$ is a discrete faithful
representation of $\pi_1(N)$, then for all $g$ in $\pi_1(M)$,
$\rho(g)$ is not elliptic and so $\tr(\rho(g)) \neq 0$.  So a
component of $X'$ containing a discrete faithful representation
satisfies $p = 0$ (the part of $X'$ where $\a = \b = 0$ turns out to be
the character variety of $N_{1/2}$).

Let $\rho$ be a representation such that $p = 0$ and $\g^2 = 4$.  Then
$\b^2 = \a^2 = 4$ and if $L$ is the image of $l = a b a^{-1} b^{-1}$,
$\tr(L) = \a^2 + \b^2 + \g^2 - \a \b \g - 2 = -2 \a^2 + 10 = 2$.  Then
$\rho$ restricted to $\pi_1(F)$ is reducible.  Recall we have only
been looking at the part of $X$ which is the closure of the subset of
characters which are irreducible when restricted to $\pi_1(F)$.  Let
$X_1$ be the closure of the part of $X$ satisfying $p = 0$ and $\a, \b
\neq 0$ and $\g^2 \neq 4$.  Henceforth, I will only be concerned with
$X_1$.  Combining $p = 0$ with $2 \a = \b \g$ we get $2\a^2 - \g^2 -
4$.  Expanding $ t a = b^{-1} a^{-1} b^{-1} t$, $b^{-1} t b = a
(b^{-1} a^{-1} b^{-1} )^3 t$, and $b t a = a^{-1} b^{-1} t = (b^{-1}
a^{-1} b^{-1} )^3 t b^{-1} = t a^3 b^{-1}$ with the identities
(\ref{idents}) we get the following linear equations for the traces
\[ \begin{array}{rrrrl}
- \b \g \tr(T) &+ 2 \tr( A T) &+\g \tr(B T)& &= 0 \\ - \a^2 \tr(T)
&+\a ( \a^2 - 2 ) \tr(A T) & & &= 0\\ \a \b (a^2 - 2) \tr(T)& &- \a (
\a^2 - 2 ) \tr(B T) &- \a^2 \tr(A B T) & = 0.
\end{array} \]
Solving this and simplifying the answer using $2 \a = \b \g$ and
$\g^2 = 2(\a^2 - 2)$ we get
\begin{equation}
\label{star}
\tr(T) \equiv \t, \tr(A T) = \frac{\b}{\g} \t , \tr(B T) = \frac{2 \b(\a^2
  - 3)}{\g^2} \t \mtext{and} \tr(A B T) = \frac{2}{\g} \t.
\end{equation}
From $t a t^{-1} b = b^{-1} a^{-1}$ we get 
\[
\tr( A B T ) \tr(T) - \tr(B T) \tr(A T) = 2 \g  - \a \b.
\]
Combining this with (\ref{star}) and using the equations relating
$\a$, $\b$, and $\g$ everything magically simplifies and we get $4
\t^2 = - \g^4$!  (Here we need $\g^2 \neq 4$.)  Thus $X_1$ breaks up
into two components, one satisfying $2 \t = i \g^2$ and the other $2
\t = -i \g^2$.  These two components are isomorphic (recall that,
fixing $A$ and $B$, $T$ is only unique up to sign).  Let $X_0$ be the
component where $2 \t = -i \g^2$.  Now, projection of $X_0$ onto the
subspace with coordinates $\a$ and $\g$ is 1-1.  Let $C$ be the curve
in this subspace defined by $2 \a^2 - \g^2 - 4 = 0$.  Since $p = 0$,
$\g \neq 0$ on $X_0$ and the image of the projection of $X_0$ is $C
\setminus \{(\pm \sqrt{2}, 0) \}$.  The closure $\bar{C}$ of $C$ in
$\CP^2$ is isomorphic to $\CP^1$ since $C$ is a conic.  Now $\bar{C}
\setminus X_0$ consists of four ideal points $(\pm \sqrt{2} : 0 : 1)$
and $(1: \pm \sqrt{2} : 0)$.  Using that
\begin{equation}\label{trlt}
\tr(L T) = \tr( A B^2 T A ) = \a \b \tr(A B T) - \a \tr(A T) - \b
\tr(B T) + \tr(T) = -4/\t,
\end{equation}
it is easy to check that the ideal points $(\pm \sqrt{2} : 0 : 1)$
have associated slope $t$, $(1: \pm \sqrt{2} : 0)$ have associated
slope $l t$, and the root of unity at each of them is $\pm i$.  Using
(\ref{trlt}) and corresponding expression for $\tr(L T^2)$, you can
check that the component of the plane curve $D_N$ of \cite{CCGLS}
corresponding to $X_0$ is defined by the following equation, where $x$
is the eigenvalue of $T$ and $y$ the eigenvalue of $L T$
\begin{equation}
x y + i x + i y + 1 = 0. \label{answer}
\end{equation}
As a double check for the computations in this section, I had
\SnapPea\ compute eigenvalues at various Dehn fillings and these
satisfied (\ref{answer}).

\noindent {\it Remark \arabic{section}.1.} A boundary slope is called
strict if it is the boundary slope of a surface which is not a fiber
of a bundle over $S^1$ or the common frontier of two twisted interval
bundles.  Boundary slopes associated to ideal points where some
peripheral element has a pole are always strict \cite[1.3.9]{CGLS}.
An interesting open question is whether all strict boundary slopes are
associated to ideal points.  Note that the boundary slope of the
surface $S_3$ in $N$ is not associated to any ideal point of a
component of $X$ which contains the character of a discrete faithful
representation.  However, the boundary slope of $S_3$ is associated to
an ideal point of another component of $X$, $X(N_{1/2})$, and hence is
strict.  Thus one can't prove that every strict boundary slope is
associated to some ideal point by looking only at the components of
$X$ which have the most geometric meaning.

\noindent {\it Remark \arabic{section}.2.} The manifold $N$ is a Dehn
filling of a two-cusped manifold $M$ of Neumann and Reid which has
strong geometric isolation.  The manifold $M$ fibers over a twice
punctured torus and is the second example in Section 3 of
\cite{NeumannReid93} (see also \cite{Calegari96}).  Dehn filling one
cusp of $M$ to get a once-punctured torus bundle yields $N$.  That $M$
has strong geometric isolation means that if you do hyperbolic Dehn
filling on one cusp and then do any hyperbolic Dehn filling on the
other cusp, you don't change the length of the geodesic that was added
to the first cusp.  Let $T_1$ and $T_2$ be the two tori boundary
components of $M$.  Let $Y_0$ be a component of the character variety
of $M$ which contains the character of a discrete faithful
representation, and $X(T_i)$ denote the character variety of $T_i$.
Consider the map $Y_0 \to X(T_1) \cross X(T_2)$ induced by the
inclusions $T_i \injects M$.  Strong geometric isolation is equivalent
to the image of this map being a product.  Consider a one-cusped
manifold $N'$ obtained by hyperbolic Dehn filling a cusp of $M$.
Since there is a self homeomorphism of $M$ which interchanges the
cusps, we can assume $N'$ is obtained by filling in the second cusp.
Since $M$ has strong geometric isolation, the image of the character
variety of $N'$, $X(N')$ in $X(T_1)$ is independent of which Dehn
filling was done to obtain $N'$.  Since $N$ is obtained in this way, a
component $X_0(N')$ of $X(N')$ which contains the character of a
discrete faithful representation is almost the same as $X_0(N)$.  In
particular, $X_0(N')$ has ideal points where the associated roots of
unity are $\pm i$.  Also, if we view $N$ as obtained by Dehn filling
the second cusp of $M$, the trace function associated to the geodesic
added to the second cusp must be constant on $X_0(N)$.  It's easy to
check from the above description of $X_0(N)$ that the trace of $a^2 t$
is always $\pm 2 i$.  It turns out that the trace of $a^2 t$ is
non-constant on $X(N_{1/2})$.  Peter Shalen points out that this begs
the following question.  Suppose $N'$ is any one-cusped hyperbolic
manifold and $\gamma$ a non-trivial element of $\pi_1(N')$.  Is there
always a component of $X(N')$ on which $\tr_\gamma$ is non-constant?
(This component of $X(N')$ should also contain the character of an
irreducible representation.)

\section{The manifold m137}\label{M}

The example with a sixth root of unity I checked rigorously is m137,
which I will call $M$. It has one cusp, and its volume is about
$3.6638$.  It is the complement of a knot in $S^2 \cross S^1$.  The
manifold $M$ is obtained by $0$ Dehn surgery on either component of
the link shown in Fig.~1(a) (this link is $7^2_1$ in
\cite{Rolfsen76}).
\begin{figure}[htb]
\centerline{\psfig{figure=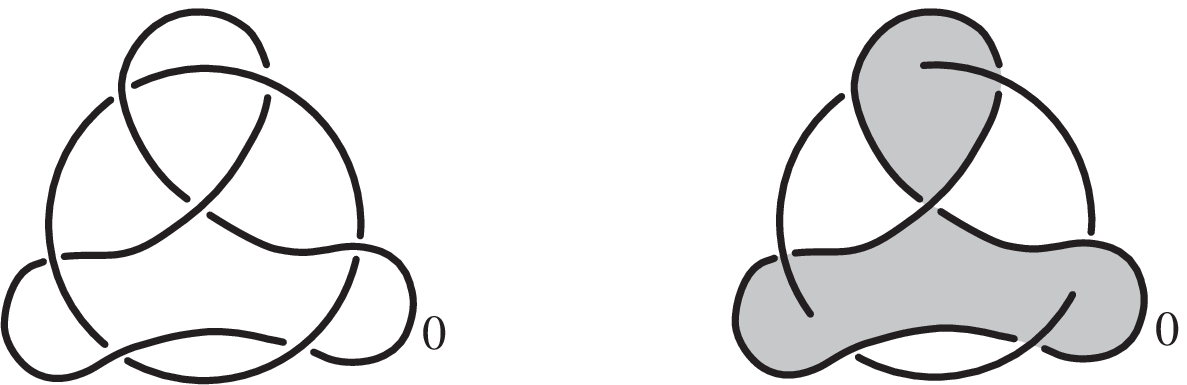, scale=1.2}}
\vskip0.1in
\hskip1.15in (a) \hskip3.3in (b)
\vskip0.1in
\noindent {\sc Fig.1.} (a) The manifold $M$ is obtained by Dehn
surgery on this link in $S^3$.  (b) The union of the pictured thrice
punctured disk with a meridian disk in the added solid torus form an
incompressible surface $F$ in $M$.  The surface $F$ is the
intersection of $M \subset S^2 \cross S^1$ with $S^2 \cross
\{\mbox{point}\}$.
\end{figure}
It consists of four ideal tetrahedra glued together as follows: Number
the tetrahedra 0, 1, 2, and 3.  Number the vertices of each
tetrahedron 0, 1, 2, and 3 in a fixed way.  We index the faces of the
tetrahedron by the opposite vertex.  The entry $k : \sigma $, where
$\sigma$ is a permutation of $\{ 0, 1, 2, 3 \}$ in the $i^{\mbox{th}}$
row and $j^{\mbox{th}}$ column of the table below means: Face $j$ of
tetrahedron $i$ is glued to face $\sigma(j)$ of tetrahedron $k$ so
that for $l \neq j$, vertex $l$ of tetrahedron $i$ is glued to vertex
$\sigma(l)$ of tetrahedron $k$ (this is the scheme used by
\SnapPea) .  
\begin{center}
\begin{tabular}{|c|c|c|c|} \hline
    1 : (0132)  & 1 : (3201) & 2 : (0132) & 3 : (0132) \\ \hline
    0 : (0132) &  2 : (1230) & 0 : (2310)  & 3 : (3201) \\ \hline
    3 : (3120) & 3 : (3012) & 1 : (3012) & 0 : (0132) \\ \hline
    2 : (1230) & 1 : (2310) & 0 : (0132) & 2 : (3120) \\ \hline
\end{tabular}
\end{center}

\SnapPea\ computes that:
\[
\pi_1(M) = \left \langle a, b \  \Big|  \ a^3 b^2 a^{-1} b^{-3} a^{-1} b^2 \right\rangle.
\] 
A basis of $\pi_1(\bdry M)$ is $\alpha = a^{-1} b^{-1}$ and $\beta =
a^{-1} b^2 a^4 b^2$ where $\alpha$ is a meridian of the knot.   I will
show there is an ideal point where the slope of an associated surface
is $\alpha$ and the associated root of unity is a sixth root not $\pm
1$.

Figure 1(b) describes an incompressible thrice punctured sphere $F$ in
$M$.  Because $\pi_1(F)$ is generated by two conjugates of $\alpha$
whose product is also a conjugate of $\alpha$, any element of
$\pi_1(F)$ has finite trace at the ideal point.  Thus, $F$ stabilizes
a vertex of the associated tree.  Presumably, the surface associated
to this ideal point is some number of parallel copies of $F$, perhaps
tubed together along $\bdry M$.

The gluing equations of $M$ are:
\begin{equation} \begin{array}{r c l}
z_1 ( 1 - z_2) z_3 (1 - z_4) &=& (1 - z_1) (1 - z_3) \\
( 1 - z_1)^2 (1 - z_2) &=& - z_1^2 z_2 (1 - z_3) (1 - z_4) \\
z_2 (1 - z_3) &=& (1 - z_1) z_4 \\
z_1 (1 - z_3) z_4 &=& -(1 - z_2)^2 z_3, \end{array} \label{glueingequations} 
\end{equation}
where $z_i$ is the edge parameter of the $i^{\mbox{th}}$ tetrahedron
for a suitable choice of edges.  This defines the deformation variety
$T \subset \C^4$.  The holonomy functions are given by:
\[
h_\alpha = - \frac{ (1 - z_2)( 1 - z_4)}{(1 - z_1) z_4 }, h_\beta =
\frac{z_1 (1 - z_4)}{(1 - z_1) z_2 (1 - z_3)^2 z_4}.
\]
The point $p_0 = (1/2 + i/2, 1 + i, 1/2 + i/2, 1 + i)$ defines a
complete hyperbolic structure on $M$. The point $p = (\zeta, 1, 1, 0
)$, where $\zeta = (1/\sqrt{3}) e^{\pi i /6}$, is our desired ideal
point.  Using the second and fourth gluing equations, we see that on
$T$ \[ h_\alpha = - \frac{(1 - z_1)}{z_1 z_2 z_3}.
\] 
Hence $h_\alpha(p) = e^{2 \pi i/3}$.  From the formula for $h_\beta$,
we see that the numerator converges to $\zeta$ as we head toward $p$, and
the denominator goes to zero, so $h_\beta$ has a pole at $p$.

It remains to show there is a sequence of points $\{ p_i \}$ of $T$
converging to $p$ where no coordinate of any $p_i$ is equal to $0,1,$
or $\infty$, since then $\{ p_i \}$ gives rise to a sequence of
representations $\rho_i \from \pi_1(M) \to \PSL(2, \C)$.  As $H_1(M) =
\Z$, every representation of $\pi_1(M)$ into $\PSL(2, \C)$ lifts to
$\SL(2,\C)$ (see the lemma in section 6 of \cite{CCGLS}), and so we have a
sequence of representations $\widetilde\rho_i$ into $\SL(2, \C)$ with
$\tr(\widetilde\rho_i(\beta) ) \to \infty$.  Passing to a subsequence,
the sequence of characters of the $\widetilde\rho_i$ converges in $X$
to an ideal point, and $\tr(\widetilde\rho_i(\beta) ) \to \lambda +
1/\lambda$ where $\lambda$ is a sixth root of unity not $\pm 1$.  This
will show that $M$ has an ideal point whose root of unity is a sixth
root not $\pm 1$.

Note that any point of the form $(z, 1, 1, 0)$ is in $T$, and any
point of $T$ near $p$ with a coordinate equal to $0,1,$ or $\infty$ is
of this form.  Producing the sequence $\{ p_i \}$ is equivalent to
showing that some other component of $T$ besides the one defined by $\{ z_2
= 1, z_3 = 1, z_4 = 0 \}$ contains $p$.  To do this we check that the
tangent space to $T$ at $p$ is two-dimensional.  It is necessary to
find the radical of the ideal defined by (\ref{glueingequations}) for
this.  I did this with Macaulay2 \cite{Macaulay2}, and the
tangent space is the same as the two-dimensional subspace which you
would get if you assumed that (\ref{glueingequations}) defines a
radical ideal.

It's interesting to note that the volume at the ideal point $p$ is
positive, since exactly one tetrahedron is not collapsed.  In most of
the examples I looked at, as you approached the ideal point all of the
tetrahedra collapsed completely and the volume went to $0$.
 

\end{document}